\documentclass[12pt]{article}
\usepackage[all]{xy}
\usepackage{graphicx}

\usepackage[latin1]{inputenc}
\usepackage[T1]{fontenc}
\reversemarginpar

\usepackage{amsmath}
\usepackage{amssymb}
\usepackage{latexsym}
\usepackage{amsfonts}
\usepackage{epsfig}
\usepackage{picinpar}
\usepackage[osf,sc]{mathpazo}
\usepackage{microtype}
\usepackage{textcomp}
\newcommand{\C}{{\mathbb C}}

\newcommand{\Q}{{\mathbb Q}}

\newcommand{\Z}{{\mathbb Z}}
\title{An e-lliptic curve}

\begin{document}
\title{A remarkable eelliptic curve}
\author{Duco van Straten \footnote{Johannes Gutenberg University, Institut f\"ur Mathematik, Staudingerweg 9, 55128 Mainz, Germany}}
\maketitle
\begin{abstract}
  We describe a system of plane algebraic curves defined over $\Z$,
  attached naturally to the exponential function. One of these is a
  remarkable curve of degree 6 that has genus equal to $1$, and looks like
  \begin{center}
    \includegraphics[width=4cm]{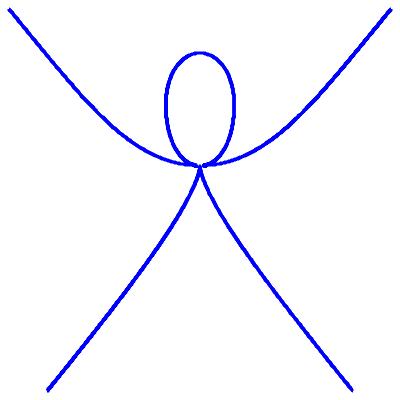}
\end{center}
As the sextic curve has rational points, it is an elliptic curve and
can be transformed over $\Q$ into the curve ${\bf 1584.j1}$ of the {\sc LMFDB}.
One is left to wonder what the number $11$, appearing in the factorisation
$1584=2^4 \cdot 3^2\cdot 11$, has to do with the exponential function.
\end{abstract}

{\bf 1.} The exponential function $e^x$ is the first object we encounter
in elementary analysis that transcends the algebraic geometrical
world of algebraic functions. Being its own derivative and having
no zeros or poles, it seems to lack clear geometrical features.
A little later we learn to extend $e^x$ into the complex domain and study the
function $e^z$, $z=x+iy$, and discover that it has an imaginary period $2\pi i$,
and is related, via
\[e^{x+iy}= e^x\left( \cos y+i\sin y \right),\]
to the functions $\sin$ and $\cos$, the circle and the rest of trigonometry and cyclotomy.\\

{\bf 2.} In the real domain, the function stretches between the limit value $0$
for $x=-\infty$ and $+\infty$ for $x=+\infty$, thus showing that
the point $z=\infty$ is an {\em essential singularity}. To study this special point,
it is convenient to perform the inversion $ z \mapsto 1/z$ and study $e^{1/z}$ in a neigbourhood of
$z=0$. Decomposing $z=x+iy$ in real and imaginary parts,
the reciprocal transformation takes the form  \[(x,y) \mapsto (x/(x^2+y^2),-y/(x^2+y^2)),\]
The level sets of $|e^z|=e^x$ are the vertical lines $x=\text{constant}$ and these
are mapped to the system of circles tangent to line $x=0$ and with center on the line
$y=0$, the level sets of  $|e^{1/z}|=e^{x/(x^2+y^2)}$.
\begin{center}
    \includegraphics[width=4cm]{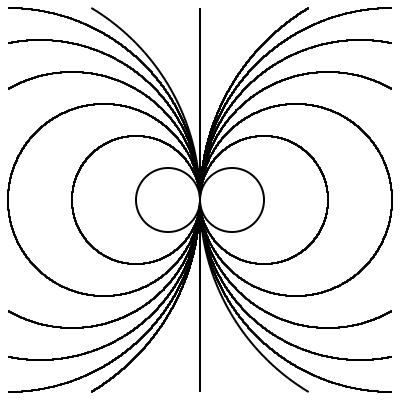}
\end{center}
If we approach the origin $z=0$ from right halfplane $x=\text{Re}(z)>0$, the function $e^{1/z}$ is unbounded; coming from the left halfplane $\text{Re}(z) < 0$, the function approches $0$. The restriction to the imaginary axis $x=0$ is $e^{1/(iy)}$ and has $\cos (1/y)$ as
real part and exhibits the well-known wild oscillatory behaviour near $0$.\\

{\bf 3.} To get a clearer view of what is happening, one can plot, for small $x=\text{constant}$,
the function $\text{Re}(e^{1/(x+iy)})=e^{x/(x^2+y^2)} \cos \frac{y}{x^2+y^2}$. For $y$ large this
function tends rapidly to its limiting value $1$ so we look at a small interval around $y=0$.
For $x<0$ one observes an oscillatory behaviour, whose amplitude suddenly drop to very small values, forming a {\em valley} around the origin. For $x>0$ we have complementary behaviour:
function is small well away from the origin, but develops oscillations
that suddenly blow up when we approach the origin.
\begin{center}
\includegraphics[width=6cm]{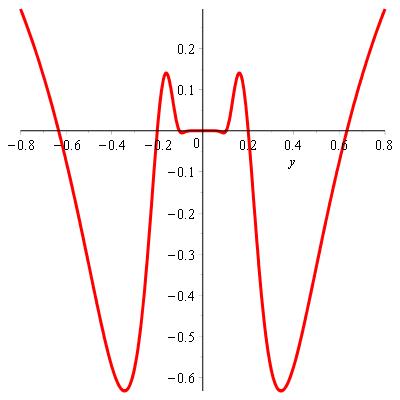} \hspace{1cm}
\includegraphics[width=6cm]{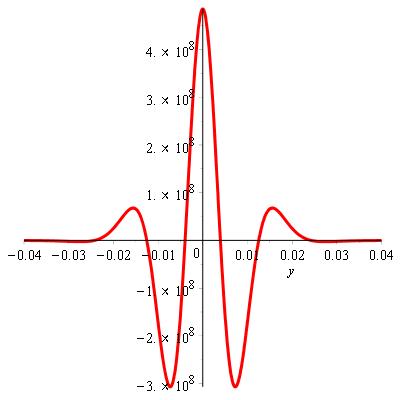}

Plot of $\text{Re}(e^{1/(x+iy)})$ for $x=-0.05$ (left) and $x=0.05$ (right).\\
(Note the difference in scaling on the axes.)
\end{center}
The graphs of $|e^{1/(x+iy)})|=e^{x/(x^2+y^2)}$ for $x=\text{constant}$ suggests
the presence of two inflectional points that may serve to demarcate the location
of the regions of there the function is very small.
\begin{center}
\includegraphics[width=6cm]{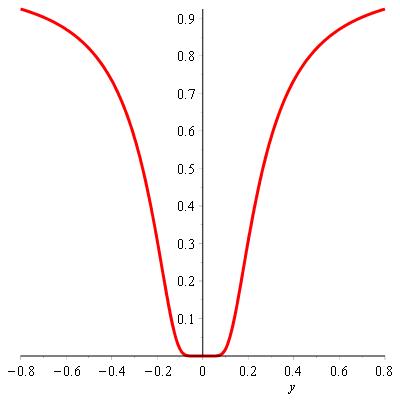} \hspace{1cm}
\includegraphics[width=6cm]{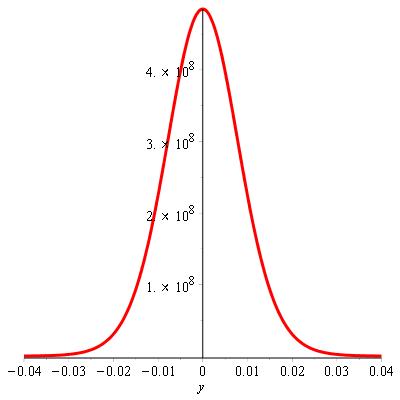}

Plot of $|(e^{1/(x+iy)}|$ for $x=-0.05$ (left) and $x=0.05$ (right).\\
(Note the difference in scaling on the axes.)
\end{center}
These pictures ilustrate the striking {\em difference in scale} for positive and
coresponding negative values of $x$: the width of the valley for $x=-c<0$ seems to
be much bigger than the width of the corresponding peak for $x=c>0$.

{\bf 4.} It is easy to calculate the location of these inflectional points by
computing the second derivative of $e^{x/(x^2+y^2)}$ with respect to $y$. Its vanishing is determined by
the polynomial factor $xF_2(x,y)$, where
\[F_2(x,y):=(x^2-3y^2)(x^2+y^2)-2xy^2,\]
which defines an irreducible quartic curve $C_2$ with a surprising
singularity of type $D_5$ at the origin, consisting of a smooth branch and a transverse
cuspidal branch. 
\begin{center}
  \includegraphics[width=6cm]{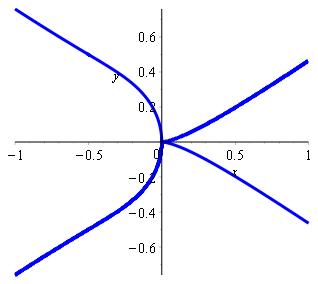}

  The rational curve $C_2$ defined by $F_2=0$.
\end{center}
Hence we are dealing with a rational curve and it can be parametrised by the pencil of lines through $0$, which intersect the curve $C_2$ in a unique further point; the substitution $y=mx$
leads to
\[x=\frac{2m^2}{(1+m^2)(1-3m^2)},\;\;\; y=m\cdot x .\]

In the halfspace $x>0$,  $C_2$ exhibits it cuspidal branch, which is intersected by
the line $x=c$ in two points, which span an interval of length
\[R_c=2\sqrt{\frac{1}{2}} c^{3/2}(1+\frac{c}{2}+\ldots).\]
whereas in the halfspace $x<0$ the curve $C_2$ is approximated near $0$ by a
parabola, and now the two points of intersection of $C_2$ with the line $x=c$
span an interval of length
$L_c =2\sqrt{\frac{2}{3}}(-c)^{1/2}(1+\frac{c}{2}+\ldots)$. The ratio between these intervals is
\[R_{c}/L_{-c} \sim \sqrt{\frac{3}{4}}\cdot c.\]
As a consequence, the central peak is indeed much narrower for small $x=c>0$ than the corresponding
central valley for $x=-c<0$, explaining the feature observed in {\bf 3.}\\

{\bf 5.} Surprised by the beauty of the curve $C_2$, one asks: why not take a look at higher derivatives of $e^{x/(x^2+y^2)}$? The third derivative factors out the irreducible polynomial
\[ F_3(x,y)= 6(x^2-y^2)(x^2+y^2)^2 +3x^5-6x^3y^2-9xy^4-2x^2y^2 .\]
The curve defined by $F_3=0$ is a sextic $C_3$, with a fourfold point at the origin, where
we find two smooth and a cuspidal branch. One of the smooth branches closes up at the point $(x,y)=(-1/2,0)$ to a loop lying in the left halfplane. 
\begin{center}
  \includegraphics[width=7cm]{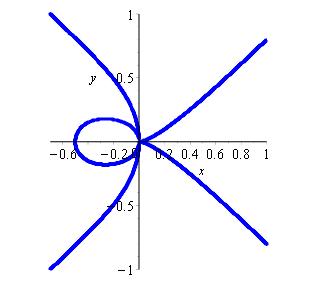}

  The curve $C_3$ defined by $F_3=0$.
\end{center}
In the affine $x,y$-plane, the curve has the origin as unique singular point.
But we will look at the projective closure of this curve, lying in the
projective closure of the complexification $\mathbb{C}^2$ of
$\mathbb{C}$ considered as $\mathbb{R}$-vector space.
As is visible from the equation, it passes doubly through the
two circle points $(\pm i: 1:0)$ on the line at infinity. As we have two
smooth and one cuspidal branch at the origin, and each pair of branches
have mutual intersection multiplicity $=2$, we obtain for the
$\delta$-invariant of $C_3$ at $0$: \[\delta(C_3,0)=1+0+0+2+2+2=7. \]
Together with the double points of $C_3$ at the circle points, we find by the
{\em Riemann-Clebsch formula:}
\[ \text{genus}(C_3)=\frac{(6-1)(5-2)}{2}-7-1-1=1,\]
so $C_3$ is  a curve of genus $1$. As we have the rational point $(-1/2,0)$
lying on $C_3$, the curve is in fact an {\em elliptic curve defined over $\Q$}.
We note the presence of the further rational points $(-1/6,\pm 1/6)$ on $C_3$.

{\bf 6.} The pencil of lines through the origin intersect the curve $C_3$
in two further points, representing $C_3$ as double cover of the projective line of
directions at $0$. Setting $y=m x$ in $F_3(x,y)$, we obtain
\[F_3(x,mx)=-x^4(Ax^2+Bx+C),\]
where 
\[A=6 (m - 1) (m + 1) (m^2  + 1)^2,\;\;\; B=3 (3 m^2 -1) (m^2 + 1),\;\;C=2m^2.\]
The $x$-coordinate of the two further intersection points of $C_3$ with the line $y=mx$
are the roots of $Ax^2+Bx+C$. So the ramification points of the double cover are
determined by the discriminant 
\[D=B^2-4AC=3(m^2+1)^2(11m^4-2m^2+3)\]
of the quadratic term. This shows that the curve $C_3$ is a double cover ramified over
the four roots of the quartic
\[ 11m^4-2m^2+3 .\]
From this we readily compute the $j$-invariant as
\[ j=\frac{62500}{33}=\frac{2^2\cdot 5^6}{3\cdot 11},\]
which determines the elliptic curve as a complex curve. However, the above
quartic has only complex roots, and it is not directly clear how to find the familiar Weierstrass normal form from here.\\

{\bf 7.} It is much more convenient to first perform the inversion $z \mapsto 1/z$, or
$(x,y) \mapsto (x/(x^2+y^2),-y/(x^2+y^2))$ on our curve, so in a way 'undoing'
the initial inversion to go from $e^z$ to $e^{1/z}$. The corresponding
substitution into $F_3$ leads to:
\[F_3\left(\frac{x}{x^2+y^2},-\frac{y}{x^2+y^2}\right)=\frac{G_3(x,y)}{(x^2+y^2)^4},\]
where
\[G_3:=-2x^2y^2+3x^3-9xy^2+6x^2-6y^2 .\]
As the inversion is an involutory transformation, we also have
\[ G_3\left(\frac{x}{x^2+y^2},-\frac{y}{x^2+y^2}\right)=\frac{F_3(x,y)}{(x^2+y^2)^4} .\]
Indeed, from the standpoint of projective geometry, we are performing a {\em Cremona
  transformation} based at $0$ and the two circle points. Any quadric through
these three points (i.e. any linear combination of $xz, yz, x^2+y^2$)
intersects $C_3$ generally in $2\cdot 6-4-2-2=4$ further points, so this
transforms our sextic $C_3$ to a {\em quartic curve $Q$}, defined by vanishing of
the above polynomial $G_3$.
\begin{center}
   \includegraphics[width=6cm]{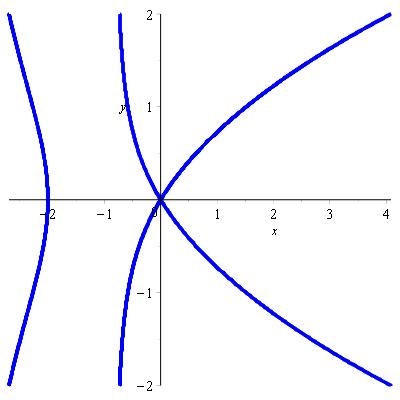}

The quartic curve $Q$ defined by $G_3=0$.
 \end{center}
The (projective closure of) the curve $Q$ has $(0:0:1)$ and $(0:1:0)$ as
singular points, so the genus is indeed $3-1-1=1$.
We note that the equation for $Q$ is of a very special form:
\[0=-(2x^2+9x+6)y^2+3x^2(x+2),\;\;\; y^2=3x^2\frac{x+2}{2x^2+9x+6} .\]

The introduction of the variable $y':=y(2x^2+9x+6)/x$ converts this into
the form
\[ y'^2= 3(x+2)(2x^2+9x+6)=6x^3+39x^2+72x+36.\]
Completing the cube and scaling $u=6x+13$, $v=6y'$ brings the equation in
{\em minimal Weierstrass form}
\[ v^2=u^3-75u+ 74,\]
defining a standard elliptic curve $E$ in the $u,v$-plane.
\begin{center}
  \includegraphics[width=6cm]{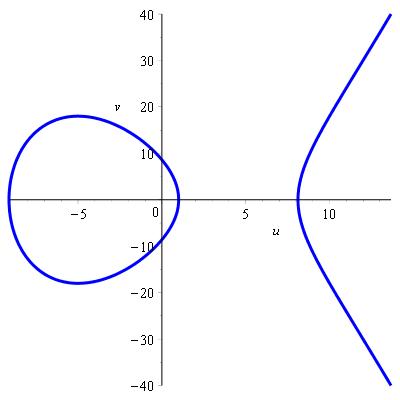}
  
The elliptic curve $E$ defined by $v^2=u^3-75u+74$.
 \end{center}

Combining the transformations, we can express $x$ and $y$ in $u$ and $v$:
\[ x=\frac{u-13}{6},\;\;\;y=\left(\frac{v}{2}\right)\frac{u-13}{u^2+u-74},\]
and indeed
\[G_3(\frac{u-13}{6},\left(\frac{v}{2}\right) \frac{u-13}{u^2+u-74})=\frac{(u-13)^2(u^3-75u+74-v^2)}{72(u^2+u-74)},\]
checking the $E$ is mapped to $Q$. The composition \[(u,v) \mapsto (x,y) \mapsto (x/(x^2+y^2),-y/(x^2+y^2))\] is the sought
for birational map from $E$ to the sextic $C_3$.\\

{\bf 7.} There are powerful software tools to handle algebraic curves and in
particular elliptic curves.

{\sc Maple} \cite{Maple} has a nice package to deal with algebraic curves, their singularities,
Puiseux expansions, their genera, etc. It has a built-in routine {\tt Weierstrassform} that
can be applied to find the Weierstrass normal form for the elliptic curve defined by the sextic $C_3$. 
However, the algorithm produces a rather complicated transformation, whose
geometric origin is obscure. The above simple algebra-geometric approach leads to a much simpler transformation.

{\sc Pari/GP} \cite{PariGP} is a powerful software, pitched toward the computation of 
arithmetical invariants of (elliptic) curves defined over number fields. It can count points, handle modular forms and $L$-functions and much more. Convenient is the command {\tt ellidentify}, which tells us
that our curve $E$ has the Cremona label $\bf 1584f1$ \cite{Cremona},
so the conductor of our curve is $1584=2^4 \cdot 3^2\cdot 11$.

The {\sc LMFDB} (L-functions and Modular Forms Data Base) \cite{LMFDB} represents a huge extension of the earlier tables by J. Cremona \cite{Cremona}; it assigns to our curve the label $\bf 1584.j1$ and presents further data attached to the curve in an attractive form. We learn that the Mordell-Weil group of the curve
is $\mathbb{Z} \oplus \mathbb{Z}/2$. The torsion point $(u,v)=(1,0)$ of $E$ corresponds
to $(-2,0)$ on $Q$ and $(-1/2,0)$ on $C_3$. The generator of infinite order $(u,v)=(-5,18)$ is mapped to $(-3,3)$ on $Q$ and to $(-1/6,1/6)$ of $C_3$. Furthermore, the curve $E$ posseses $11$ (!) points with integral coordinates, which map to special rational points on $Q$ and $C_3$.
From {\sc LMFDB} we also learn that our curve has the curve $\bf 264.c1$ with eqation
$y^2=x^3+x^2-8x$ as minimal quadratic twist and that this curve has Mordell-Weil group $\mathbb{Z}/2$.

For general background information on algebraic curves, among the many books available,
 one of our favourites is \cite{BK}. For the arithmetic theory of elliptic curves \cite{Silverman} is indispensable.\\

{\bf 8.} Why stop at the third derivative?
Write $\phi:=\exp(x/(x^2+y^2))$ and $\partial:=\partial/\partial y$, so
\[ \partial \phi= h\phi,\;\;\; h:=\partial (\frac{x}{x^2+y^2}).\]
For the $n$-th derivative one then has
\[ \partial^n \phi=h_n \phi,\;\;\;\;h_{n+1} =\partial h_n +h \cdot h_n,\]
so that
\[ h_n=\frac{f_n}{(x^2+y^2)^{2n}},\]
where $f_n\in \Z[x,y]$ is a polynomial of degree $3n-1$. For $n$ odd $f_n$
factorises as  $f_n=xyF_n$, and for $n$ even as $f_n=xF_n$, where $F_n \in \Z[x,y]$ is an irreducible polynomial (of degree $3n-3$ (n odd) resp. $3n-2$ (n even)).(Note $f_0=1$, $f_1=-2xy$, so there are no curves $C_0$ or $C_1$.)

\begin{center}
  \includegraphics[width=6cm]{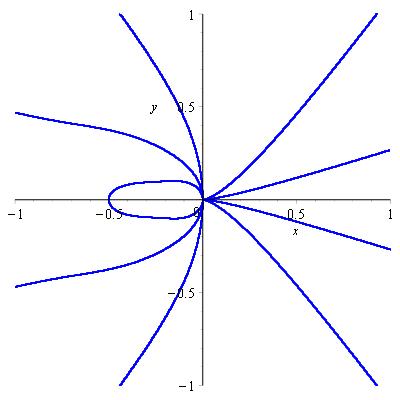} \hspace{2cm}
  \includegraphics[width=6cm]{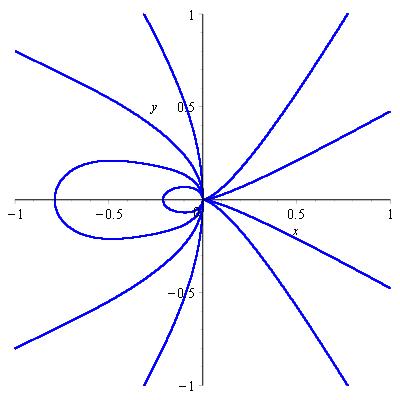}
  
The curve $C_4$ of genus $4$ and $C_5$ of genus $7$.  
\end{center}
The curve $C_n$ defined by $F_n=0$ has a unique singularity at the origin,
where it has $(n-1)$ smooth local branches with common vertical tangent and
$\lfloor \frac{n}{2} \rfloor$ cuspidal branches with common horizontal tangent.
Furthermore it has ordinary $(n-1)$-fold points at the circular points.
The genus of $C_n$ can be computed to be
\[\text{genus}(C_n) = \left\{\begin{array}{l} (n-2)(3n-4)/4\;\;\textup{if $n$ is even}\\
 (3n^2+1)/4\;\;\textup{if $n$ is odd}\\\end{array}\right. \]

In the real domain, the curves $C_n$  have an insect-like appearance.
The $(n-1)$ smooth branches are all in the halfspace $x \le 0$;
$\lfloor \frac{n-1}{2}\rfloor$ close up to 'heads', whereas the remaining branches
form pairs of 'arms'; the cuspidal branches are in the halfspace $x \ge 0$ and
form the 'legs'.\\

{\bf 9.}
Of course, we do not need to restrict $y$ derivatives only. It is simple to show
that 
\[ P_{a,b}:=(x^2+y^2)^{2(a+b)} \left( \frac{\partial^{a+b}}{\partial x^a y^b} e^{x/(x^2+y^2)}\right)e^{-x/(x^2+y^2)}\]
is a polynomial of degree $3(a+b)-1$. For $b$ an odd number, $P_{a,b}$ is divisible by $y$; for $a=0$ there is an additional factor $x$. Clearing out these factors, we obtain a polynomial $F_{a,b}$ that defines an irreducible curve $C_{a,b}$ and $C_n=C_{0,n}$.
The curve $C_{a,b}$ has the origin as only singular point in the affine plane and
furthermore ordinary $a+b$-fold singularities at the two circle points. The genus is basically a quadratic function of $a$ and $b$, to be more precise we have
\[
\text{genus}(C_{a,b})
= 9a^2/8 + 3ab + 3b^2/4 - 13a/4 - 5b/2 +
\begin{cases}
2 & \text{if}\; a\equiv0\bmod2, \\
13/8 & \text{if}\; a\equiv0\bmod2, \; b\equiv0\bmod4, \\
9/8 & \text{if}\; a\equiv0\bmod2, \; b\equiv2\bmod4,
\end{cases}
\]
when $b$ is even and
\[
\text{genus}(C_{a,b})
= a^2 + 3ab + 3b^2/4 - 7a/2 - 3b +
\begin{cases}
7/4 & \text{if}\; a\equiv1\bmod2, \\
13/4 & \text{if}\; a=0, \\
9/4 & \text{if}\; a\equiv0\bmod2, \; a>0,
\end{cases}
\]
when $b$ is odd.
%(I am grateful to  W. Zudilin for providing these explicit formulas).
Here a table of genera of the curves appearing for small $a,b \le 5$.
%\begin{center}
%  \[
%  \begin{array}{|c|cccccc|}
%    \hline
%    \text{genus}&a=0&1&2&3&4&5\\
%    \hline
%    b=0&- &- &0 &2 &7 &13\\
%      1&- &0 &3 &7 &14&22\\
%      2&0 &3 &10&18&29&41\\
%      3&1 &6 &15&25&38&52\\
%      4&4 &13&26&40&57&75\\
%      5&7 &18&33&49&68&88\\
%    \hline
%  \end{array}
%  \]
%  \end{center}

\begin{center}
 \hspace{2cm} \hspace{3cm}
\includegraphics[width=2cm]{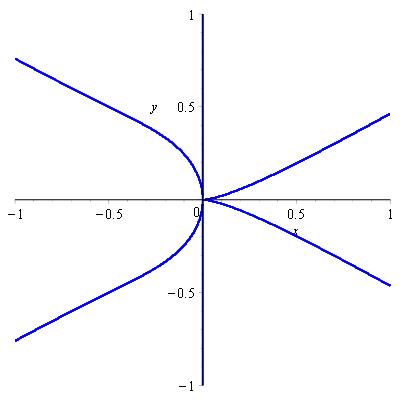}\hspace{5mm}
\includegraphics[width=2cm]{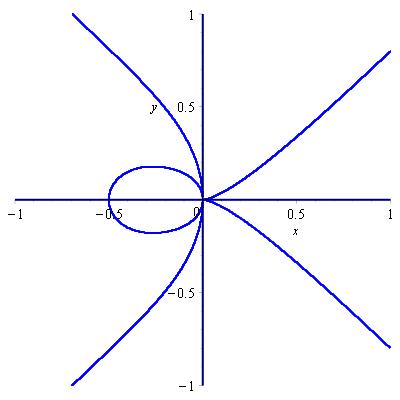}\hspace{5mm}
\includegraphics[width=2cm]{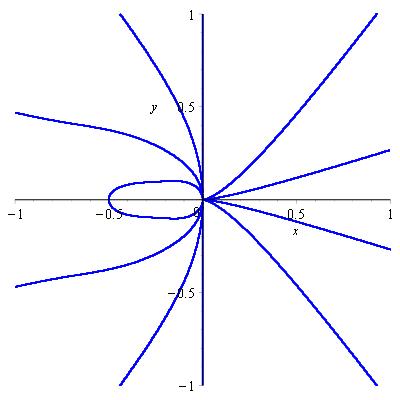}\hspace{5mm}
\includegraphics[width=2cm]{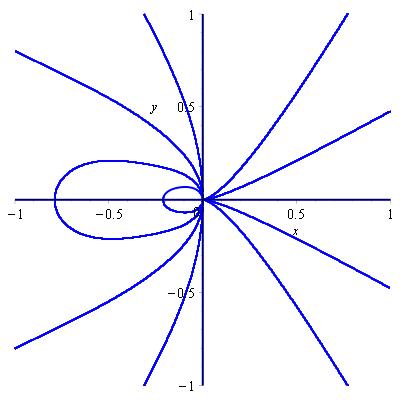}\hspace{5mm}
{\mbox{}\hfill \small\;\;\;\; $g(C_{02})=0,\hspace{5mm} g(C_{03})=1,\hspace{5mm} g(C_{04})=4,\hspace{5mm} g(C_{05})=7.$\hspace{5mm}}

\hspace{2.2cm}
\includegraphics[width=2cm]{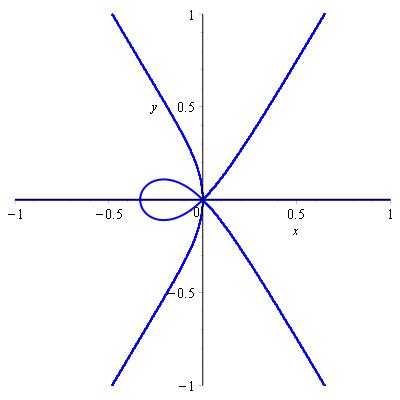}\hspace{5mm}
\includegraphics[width=2cm]{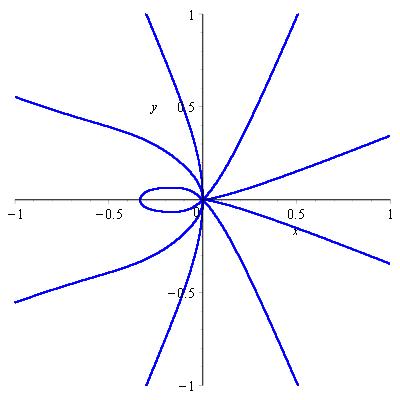}\hspace{5mm}
\includegraphics[width=2cm]{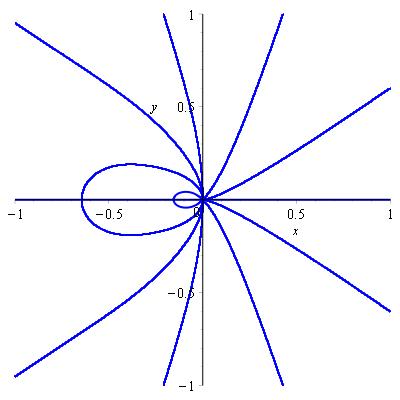}\hspace{5mm}
\includegraphics[width=2cm]{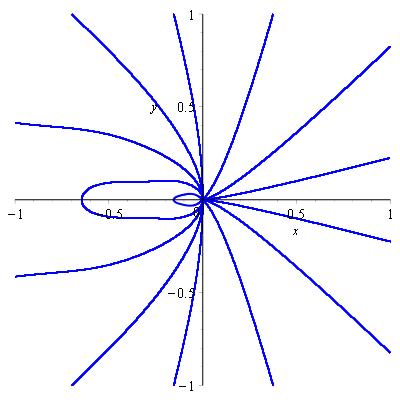}\hspace{5mm}
\includegraphics[width=2cm]{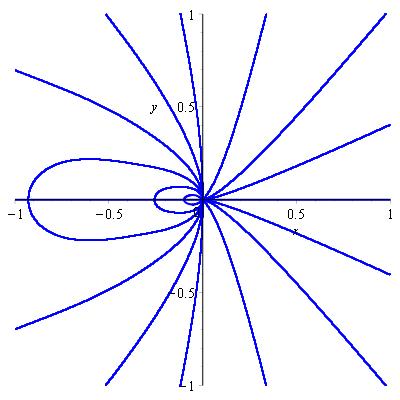}\hspace{5mm}

{\mbox{}\hfill \small\;\;\;\; $g(C_{11})=0,\hspace{5mm} g(C_{12})=3,\hspace{5mm} g(C_{13})=6,\hspace{5mm} g(C_{14})=13,\hspace{5mm}g(C_{15})=18.$ \hspace{5mm}}

\includegraphics[width=2cm]{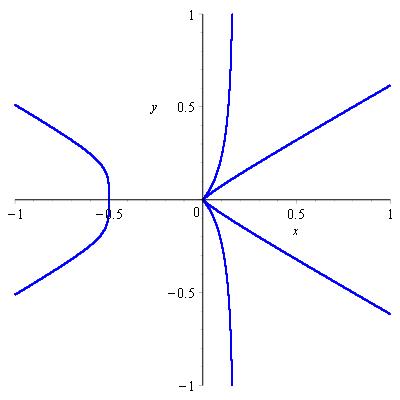}\hspace{5mm}
\includegraphics[width=2cm]{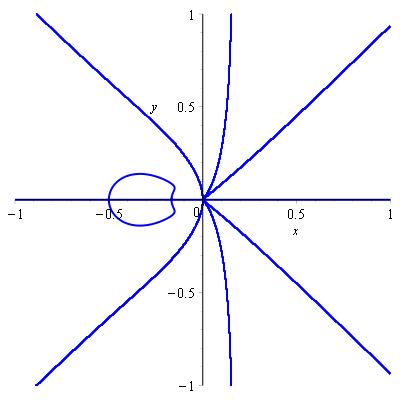}\hspace{5mm}
\includegraphics[width=2cm]{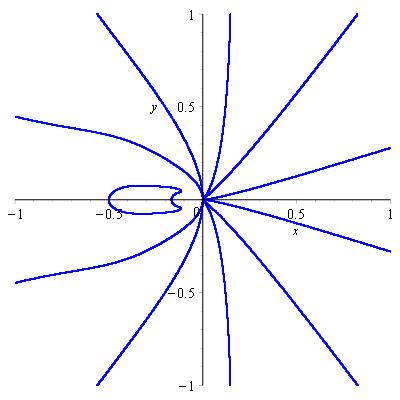}\hspace{5mm}
\includegraphics[width=2cm]{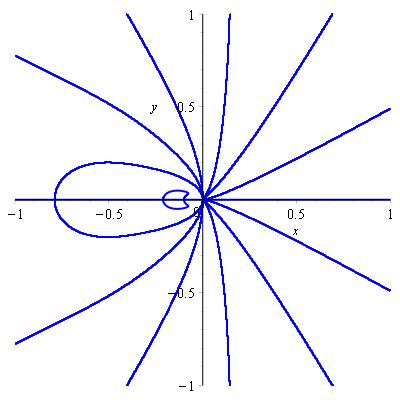}\hspace{5mm}
\includegraphics[width=2cm]{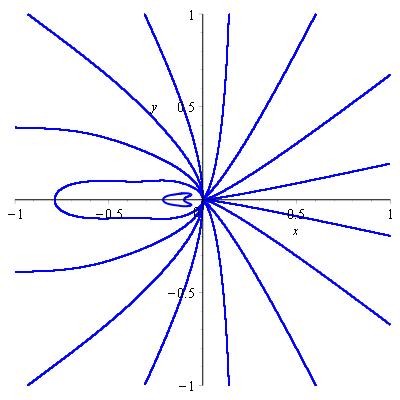}\hspace{5mm}
\includegraphics[width=2cm]{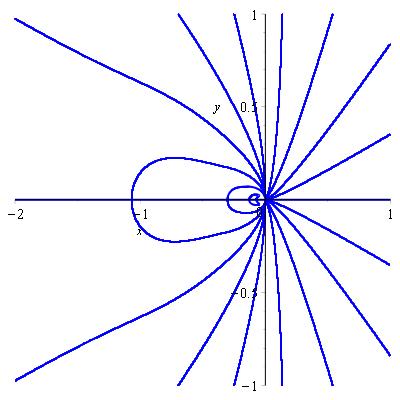}\hspace{5mm}
{\small $g(C_{20})=0,\hspace{5mm} g(C_{21})=3,\hspace{5mm} g(C_{22})=10,\hspace{5mm} g(C_{23})=15,\hspace{5mm} g(C_{24})=26,\hspace{5mm} g(C_{25})=33.$}

\includegraphics[width=2cm]{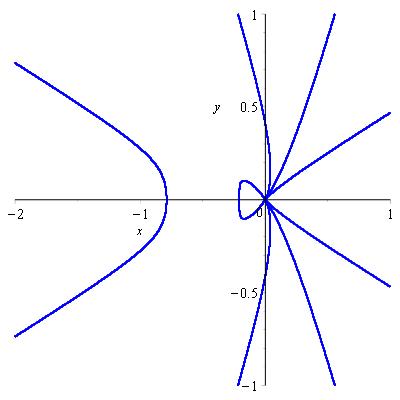}\hspace{5mm}
\includegraphics[width=2cm]{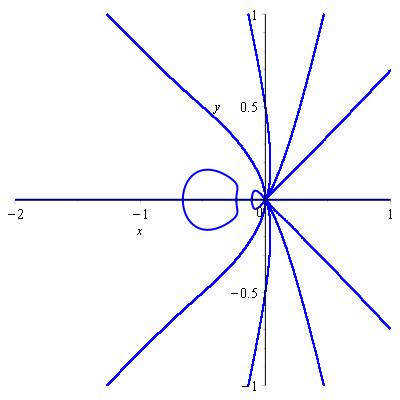}\hspace{5mm}
\includegraphics[width=2cm]{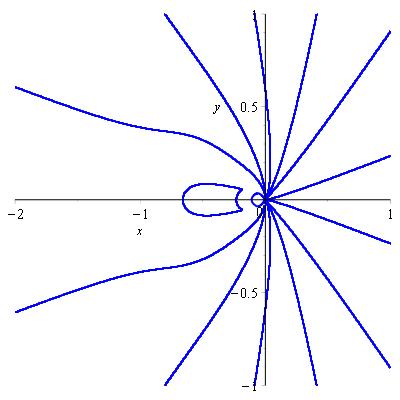}\hspace{5mm}
\includegraphics[width=2cm]{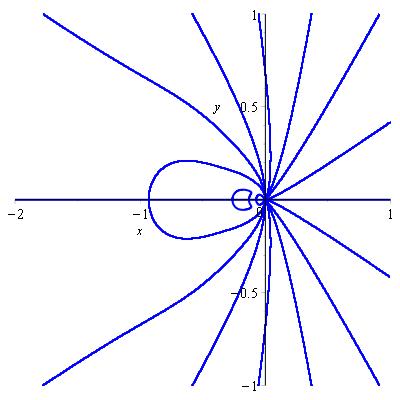}\hspace{5mm}
\includegraphics[width=2cm]{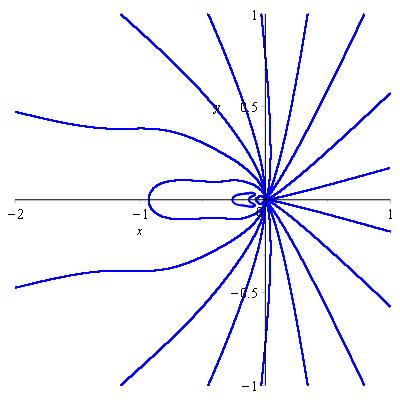}\hspace{5mm}
\includegraphics[width=2cm]{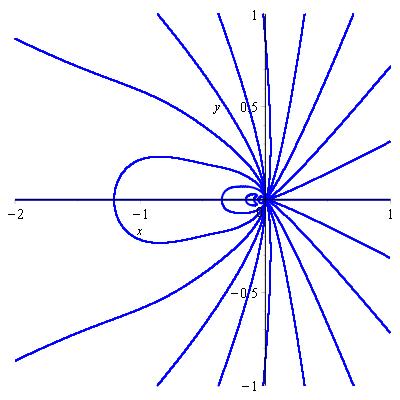}\hspace{5mm}
{\small $g(C_{30})=2,\hspace{5mm}g(C_{31})=7,\hspace{5mm} g(C_{32})=18,\hspace{5mm} g(C_{33})=25,\hspace{5mm} g(C_{34})=40,\hspace{5mm} g(C_{35})=49.$}

\includegraphics[width=2cm]{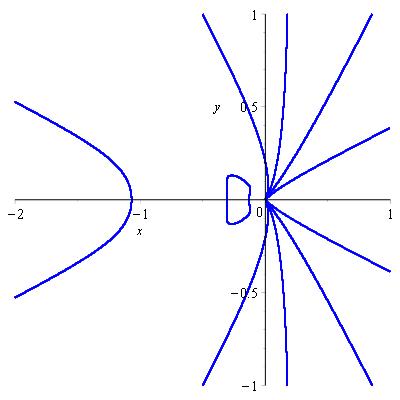}\hspace{5mm}
\includegraphics[width=2cm]{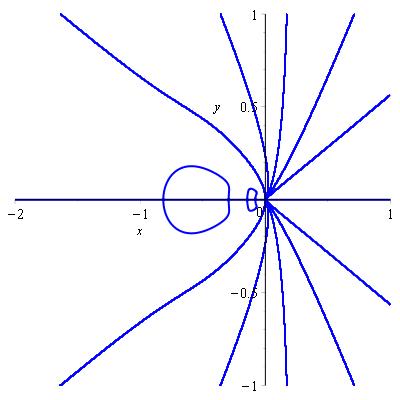}\hspace{5mm}
\includegraphics[width=2cm]{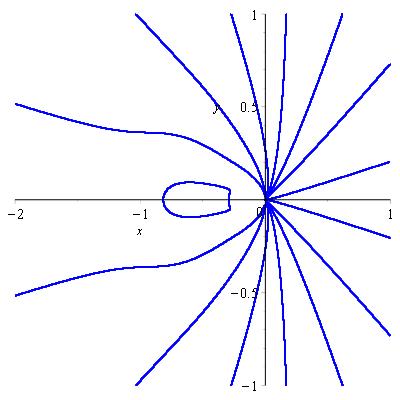}\hspace{5mm}
\includegraphics[width=2cm]{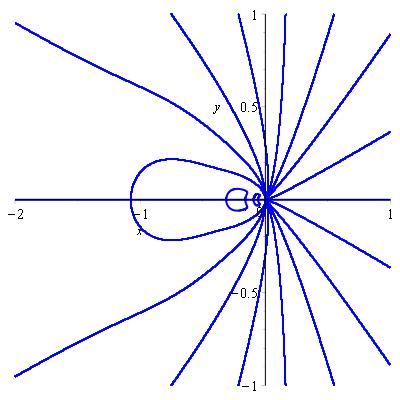}\hspace{5mm}
\includegraphics[width=2cm]{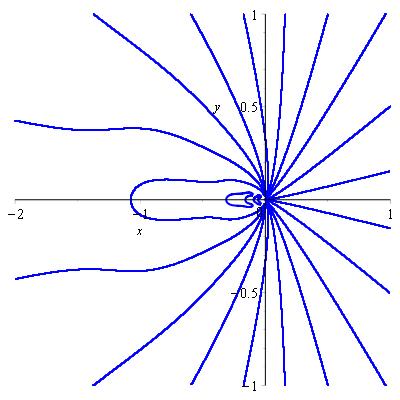}\hspace{5mm}
\includegraphics[width=2cm]{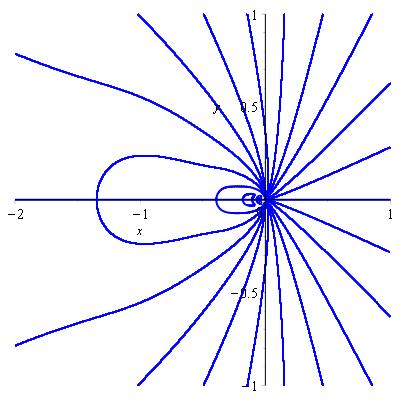}\hspace{5mm}
{\small $g(C_{40})=7,\hspace{5mm} g(C_{41})=14,\hspace{5mm} g(C_{42})=29,\hspace{5mm} g(C_{43})=38, \hspace{5mm}g(C_{44})=57,\hspace{5mm}g(C_{45})=68.$}

\includegraphics[width=2cm]{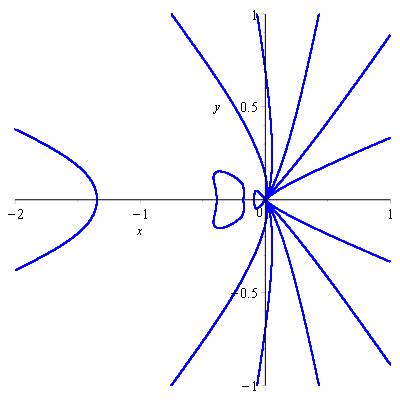}\hspace{5mm}
\includegraphics[width=2cm]{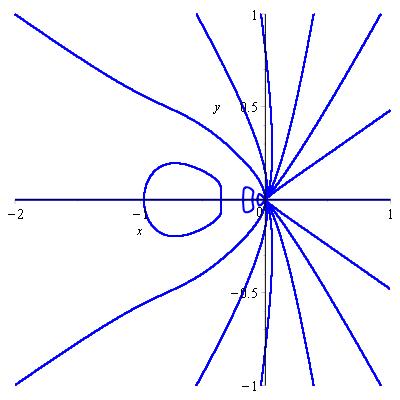}\hspace{5mm}
\includegraphics[width=2cm]{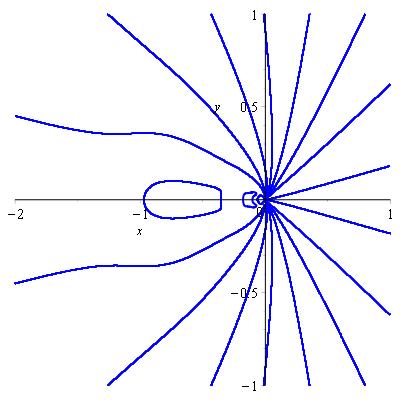}\hspace{5mm}
\includegraphics[width=2cm]{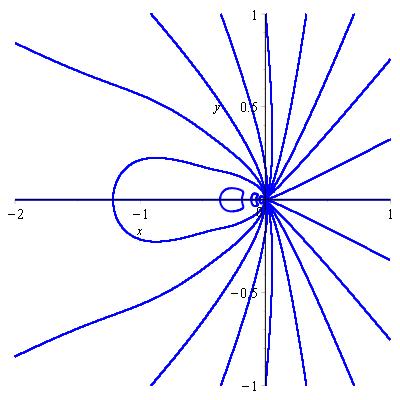}\hspace{5mm}
\includegraphics[width=2cm]{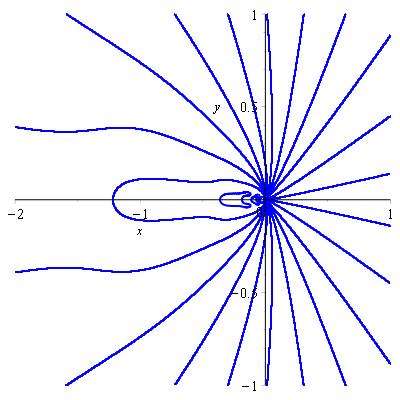}\hspace{5mm}
\includegraphics[width=2cm]{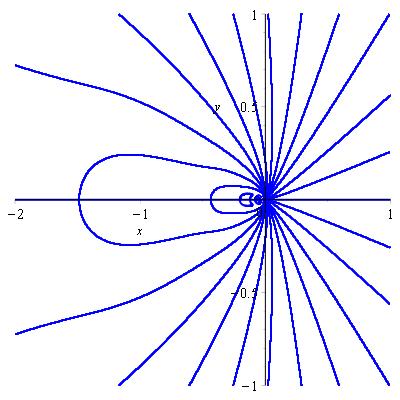}\hspace{5mm}
{\small $g(C_{50})=13,\hspace{5mm}g(C_{51})=22,\hspace{5mm} g(C_{52})=41, \hspace{5mm}g(C_{53})=52,\hspace{5mm} g(C_{54})=75,\hspace{5mm} g(C_{55})=88.$}
\end{center}

  {\bf 10.} It somewhat of a mystery to me how the curve $C_3$, with its rich specific
  arithmetic geometry, like its Galois representation of conductor $1584=2^43^211$,
  can arise without any specific choices from the exponential function, a basic function
  usually considered as belonging to ``algebraic geometry over $\mathbb{F}_1$''.
  Of course there are many other algebraic geometric objects one can 'naturally' attach to
  the exponential function. For example, we can look at the Taylor polynomial 
\[ f_n(x):=1+x+\frac{1}{2!} x^2+\ldots+\frac{1}{n!}x^n\]
of degree $n$ of $e^x$. The location of the roots of $f_n(nx)$ and their convergence
to the {\em Szeg\H{o} curve} defined by $\{z \in \C\; |\;\; |ze^{1-z}|=1\}$
is a well-known topic in classical analysis \cite{Szego}.
One can look at the sequence of hyperelliptic curves $y^2=f_n(x)$,
which for $n=3,4$ lead to two further 'natural' elliptic curves. But, in a sense
I find hard to make precise, these curves look more artificial to me.

Rather it appears that the curves $C_{a,b}$ are natural {\em algebraic geometric sattelites} of the exponential motive
defined by $e^{1/z}$, and as such has some obvious generalisations. One may look at a similar system of curves for $e^{1/z^2}$.
Here one also finds an elliptic curve, ${\bf 576.c3}$, with Mordel-Weil group $\Z \oplus \Z/2\oplus \Z/2$.  
In general, if $g, S \in K:=\Q(x_1,x_2,\ldots,x_n)$ are two rational functions, $\partial \in Der(K)$ a derivation of $K$,
and $\phi:=g. e^S$, then repeated differentiation leads a sequence $g_1,g_2,g_3,\ldots,g_n,\ldots \in \Q(x_1,x_2,\ldots,x_n)$
defined recursively by
\[g_1:=g, g_2=\partial(g_1)+g_1 \cdot \partial S,\ldots, \;g_{n+1}:=\partial g_n+g_n \cdot \partial S,\ldots\]
and which appear as prefactors of  $\partial^n\phi=g_n \phi$. The classical motives attached
to the $g_n$  are in a similar way algebraic-geometrical satellites of the exponential motive with exponential
periods $\int e^S \omega$,  $\omega:=g\cdot dx_1dx_2\ldots dx_n$, \cite{KZ}.\\
  
{\bf Acklowledgement:} I thank W. Zudilin for interest in the curve $C_3$ and many useful comments on
this note and the suggestion to write it. Furthermore,  I take the opportunity to thank the
developpers of {\sc Maple}, {\sc Pari/GP}, {\sc LMFDB} which made this work possible.

\end{document}